\title{On graphs with subgraphs of large independence numbers}
\author{Noga Alon \thanks{
Schools of Mathematics and Computer Science, Raymond and Beverly
Sackler Faculty of Exact Sciences, Tel Aviv University, Tel Aviv
69978, Israel.
Email: nogaa@tau.ac.il. Research supported in part by the Israel Science
Foundation, by a USA-Israeli BSF grant
and by the Hermann Minkowski Minerva Center for Geometry
at Tel Aviv University.} 
\and Benny Sudakov
\thanks{ Department of Mathematics, Princeton University, Princeton, NJ 08544, and
Institute for Advanced Study, Princeton. E-mail:
bsudakov@math.princeton.edu.
Research supported in part by NSF CAREER award DMS-0546523, NSF grant
DMS-0355497, USA-Israeli BSF grant, Alfred P. Sloan fellowship, and
the State of New Jersey.
} }

\date{}
\documentclass [11pt]{article}
\usepackage{graphics}
\usepackage{amsfonts}
\usepackage{amsmath}
\newtheorem{theo}{Theorem}[section]

\newtheorem{lemma}[theo]{Lemma}

\begin{document}
\maketitle
\begin{abstract}
Let $G$ be a  graph on $n$ vertices in which every
induced subgraph on $s=\log^3 n$ vertices has an independent set of size at least
$t=\log n$.
What is the largest $q=q(n)$ so that every such $G$
must contain an independent set of size at least $q$ ?
This is
one of several related questions raised by Erd\H{o}s and Hajnal. We show that
$q(n)=\Theta(\log^2 n/\log \log n)$, investigate the more general problem obtained by
changing the parameters $s$ and $t$, and discuss the connection to a related
Ramsey-type problem.
\end{abstract}

\section{Introduction}
What is the largest $f=f(n)$ so that every graph $G$ on $n$ vertices in which every
induced subgraph on $\log^2 n$ vertices has an independent set of size at least 
$\log n$, must contain an independent set of size at least $f$ ? This is one of several
related questions considered by Erd\H{o}s and Hajnal in the late 80s. The question 
appears in \cite{Er}, where Erd\H{o}s  mentions that they thought that 
$f(n)$ must be at least $n^{1/2-\epsilon}$, but they could not even prove
that it is at least $2 \log n$. As a special case of our main results here
we determine the asymptotic behavior of $f(n)$ up to  a factor of $\log \log n$,
showing that in fact it is much smaller than one may suspect (and yet
much bigger than $\log n$):
\begin{equation}
\label{e11}
\Omega(\frac{\log^2 n}{\log \log n}) \leq f(n) \leq O(\log^2 n).
\end{equation}
Another specific variant of  the above question, discussed in \cite{Er}, is
the problem of estimating the largest $q=q(n)$ so that 
every graph $G$ on $n$ vertices in which every
induced subgraph on $\log^3 n$ vertices has an independent set of size at least 
$\log n$, must contain an independent set of size at least $q$.
Here, too, one may tend to believe that $q(n)$ is large, and specifically
it is mentioned in \cite{Er} that probably $q(n)> \log^3 n$, but the correct
asymptotic behavior of $q(n)$ is smaller. In this case, our results determine the
asymptotic behavior of $q(n)$ up to a constant factor, implying that
\begin{equation}
\label{e12}
q(n)=\Theta(\frac{\log^2 n}{\log \log n}).
\end{equation}

Both problems above are special instances of the general problem of
understanding the asymptotic behavior of the function $f(n,s,t)$
defined as follows.
For $n>s>t$, let $f=f(n,s,t)$ denote the largest integer $f$ so that 
every graph $G$ on $n$ vertices in which every
induced subgraph on $s$ vertices has an independent set of size at least 
$t$, must contain an independent set of size at least $f$. 
In this note we investigate the asymptotic behavior of $f$, and obtain
rather tight bounds for this behavior for most interesting values of the
parameters. Our results provide much less satisfactory  information
about a closely related Ramsey-type 
problem of Erd\H{o}s and Hajnal discussed in \cite{Er}, which
is the following. For which functions $h(n)$ and $g(n)$, where
$n>g(n) \geq h(n)^2 \gg 1$, is there a graph on $n$ vertices
in which every induced subgraph on $g(n)$ vertices contains a clique of size $h(n)$ as
well as an independent set of size $h(n)$ ? In particular, Erd\H{o}s and Hajnal
conjectured that there is no such graph for $g(n)=\log^3 n$ and $h(n)=\log n$;
our results here do not settle this conjecture, and only
suffice to show that there is no such graph with $g(n)= c \log^3 n/\log \log n$
and $h(n)=\log n$, for some absolute positive constant $c$.

The rest of this note is organized as follows. In Section 2 we state our main results
concerning the behavior of the function $f(n,s,t)$. The proofs are described in
Section 3. The final Section 4 contains a few remarks, including the (simple) connection
between the  study of $f$ and the Ramsey-type question discussed above.

Throughout the note we omit all floor and ceiling signs, whenever these are not 
crucial. We always assume that the number $n$ of vertices of the graphs considered here
is large.  All logarithms are in the natural base $e$.
%, unless otherwise specified.

\section{The main results}

The following two theorems  provide lower bounds for the independence numbers
of graphs in which every induced subgraph of size $s$ contains an independent
set of size $t$.

\begin{theo}
\label{lower-bound1} 
Let $t< s< n/2$, 
%$t \rightarrow \infty$ 
and
let $G$ be a graph of order $n$ such that every induced subgraph
of $G$ on $s$ vertices contains an independent set of size $t$.
Denote $k=\lfloor \frac{s}{t-1}\rfloor$. Then $G$ contains an
independent set of size at least $\Omega\left( k n^{1/k} \right)$
if $k \leq 2 \log n$ and of size at least $\Omega\left(\frac{\log
n}{\log (k/\log n)}\right)$ if $k>2 \log n$.
\end{theo}

\begin{theo}
\label{lower-bound2} 
Let $2t\leq s<n/2$, 
%$t \rightarrow \infty$
and let $G$ be a graph of order $n$ such that every induced
subgraph of $G$ on $s$ vertices contains an independent set of
size $t$. Then $G$ contains an independent set of size at least
$\Omega\left( t\,\frac{\log (n/s)}{\log (s/t)}\right)$.
\end{theo}

The next result shows that there are graphs  with relatively small independence 
numbers
in which every induced subgraph of size $s$ contains an independent
set of size $t$.

\begin{theo}
\label{upper-bound1}\, For every sufficiently large $t$ and $2t
\leq s \leq n/2$ there exists a graph $G$ on $n$ vertices with
independence number 
$$ \alpha(G) \leq O\left(t\Big(\frac{n}{s}\Big)^{2t/(s-t)}\log
(n/t) \right)$$ such that every induced subgraph of $G$ of order
$s$ contains an independent set of size $t$.
\end{theo}

For certain values of $t$ and $s$ one can improve the previous
result as follows.

\begin{theo}
\label{upper-bound2}\, Let $t < s \leq n/2$, where $s \leq e^{2t}$, and 
assume, further, that
either there exists a constant $\delta >0$ such that
$(s/t)^{1-\delta} \geq \log n$ or $s/t=\Omega(\log n)$ and there
exists  a constant $\gamma >0$ such that $\log t \geq \log^{\gamma}
n$ . Then there exists a graph $G$ on $n$ vertices with
$$ \alpha(G) \leq O\left(\frac{t}{\log (s/t)}\log (n/t) \right)$$
such that every induced subgraph of $G$ of order
$s$ contains an independent set of size $t$.
\end{theo}

\section{The proofs}

\noindent 
{\bf Proof of Theorem \ref{lower-bound1}.}\, 
Suppose that $G$ contains $t-1$ disjoint
cliques whose union $U$ has size at least $s$. Then the size of
the largest independent set in $U$ is bounded by $t-1$, since
an independent set can intersect each of the $t-1$ cliques in at most one
vertex. This contradicts the property of $G$. Therefore if $G'$ is a
graph obtained from $G$ by deleting the vertices of $t-1$ disjoint
cliques with maximum union we have that $|V(G')| > n-s\geq
n/2$. We also have that the largest clique in $G'$ has size at
most $k=\lfloor \frac{s}{t-1}\rfloor$. Otherwise, by the above
discussion, $G$ will contain $t-1$ disjoint cliques each of size
at least $k+1$, whose union has at least $(k+1)(t-1)\geq s$
vertices. To finish the proof we apply the classical bound of
Erd{\H{o}}s and Szekeres \cite{ES} (see also \cite{GRS}) for the usual graph Ramsey 
numbers. This result
asserts that the maximum possible number of vertices in a graph with neither a clique of
size $k+1$ nor an independent set of size $\ell+1$ is at most
${k+\ell \choose k}$. By this estimate, $G'$ contains an
independent set of size $\ell$, where ${k+\ell \choose k} \geq
n/2$. Thus $\left(\frac{e(k+\ell)}{k}\right)^{k}>n/2$
and therefore if $k \leq 2 \log n$ we get $\alpha(G) \geq \ell
\geq\Omega\Big(kn^{1/k}\Big)$. On the other hand, if $k>2\log n$
we use the estimate $\left(\frac{e(k+\ell)}{\ell}\right)^{\ell}
\geq {k+\ell \choose \ell}= {k+\ell \choose k}\geq n/2$. This
gives that $\ell \geq \Omega\left(\frac{\log n}{\log (k/\log n)}\right)$ and
completes the proof of the theorem.
 \hfill $\Box $

\noindent 
{\bf Proof of Theorem \ref{lower-bound2}.}\, 
Let $r$ be the largest integer such that
$$n\, \left(\frac{t^2}{4e^2s^2}\right)^{r-1} \geq s.$$
It is easy to check that $r=\Omega\left(\frac{\log (n/s)}{\log
(s/t)}\right)$. To prove the theorem we construct a sequence
of pairwise disjoint independent sets $X_1,\ldots, X_r$ together with a
sequence of nested subsets $V_0=V(G)\supset V_1\supset \ldots
\supset V_{r-1}$ such that the following holds. Each $X_i$ is
a subset of $V_{i-1}$ of size $t/2$, $X_i \cap V_i =\emptyset$,
there are no edges between
$X_i$ and $V_i$ and $|V_i| \geq \frac{t^2}{4e^2s^2}|V_{i-1}|$ for
all $1 \leq i \leq r-1$. Then the union of all sets $X_i$ forms an
independent set in $G$ of size $rt/2=\Omega\left( t \frac{\log
(n/s)}{\log (s/t)}\right)$.

Assuming the sets $X_j, V_j$ have already been constructed for all
$j<i$, construct $X_i$ and $V_i$ as follows. Let $m$ be the
size of $V_{i-1}$. The inductive hypothesis implies that
$$m \geq \left(\frac{t^2}{4e^2s^2}\right)^{1-i}|V_0|
\geq n\,\left(\frac{t^2}{4e^2s^2}\right)^{r-1} \geq s.$$ Since
every subset of order $s$ in $V_{i-1}$ contains $t$ independent
vertices and every independent set of size $t$ belongs to at most
${m-s \choose s-t}$ subsets of size $s$ we have that $V_{i-1}$
contains at least ${m \choose s}/{m-t \choose s-t}$ independent
sets of size $t$. Therefore, using that $m\geq s\geq 2t$ and
${a\choose b} \leq \big(\frac{ea}{b}\big)^b$, we conclude that
there exist a subset $X_i$ of $V_{i-1}$ of size $t/2$ which is
contained in at least
\begin{eqnarray*}
 \frac{{m \choose s}}{{m-t \choose s-t} {m \choose t/2}} &=&
 \frac{m!}{(m-t)!}\frac{(s-t)!}{s!} {m \choose t/2}^{-1}
 \geq \left(\frac{m}{s}\right)^{t}
 \left(\frac{t}{2em}\right)^{t/2}\\
 &=& 
 \left(\frac{m t}{2es^2}\right)^{t/2}
 =\left(\frac{e\big(\frac{t^2}{4e^2s^2}\,m\big)}{t/2}\right)^{t/2}\\
 &\geq& {\frac{t^2}{4e^2s^2}\,m \choose t/2}
 \end{eqnarray*}
independent sets of size $t$. This implies that $V_{i-1}$ contains
at least ${\frac{t^2}{4e^2s^2}\,m \choose t/2}$ subsets of size
$t/2$ whose union with $X_i$ forms an independent set of size $t$.
Let $V_i$ be the union of all these subsets. By definition, for
every vertex of $V_i$ there is an independent set which contains
it together with $X_i$, so there are no edges from $X_i$ to $V_i$.
Also, it is easy to see that the size of $V_i$ must be at least
$\frac{t^2}{4e^2s^2}\,m=\frac{t^2}{4e^2s^2}\,|V_{i-1}|$. This
completes the construction step and the proof of the theorem.
 \hfill $\Box$

\noindent 
{\bf Proof of Theorem \ref{upper-bound1}.}\, 
We prove the theorem by considering an appropriate random
graph. As usual, let $G_{n,p}$ denote the probability space of all
labeled graphs on $n$ vertices, where every edge appears randomly
and independently with probability $p=p(n)$. We say that the
random graph possesses a graph property $\cal P$ {\em almost
surely}, or a.s.\ for brevity, if the probability that $G_{n,p}$
satisfies $\cal P$ tends to 1 as $n$ tends to infinity. Clearly,
it is enough to show that there is a value of the edge probability $p$
such that $G_{n,p}$ satisfies the assertion of the theorem with
positive probability.

Let $p=\frac{1}{4e^3t}\left(\frac{n}{s}\right)^{-2t/(s-t)}$. We
claim that almost surely every subset of $G=G_{n,p}$ of size $s$
spans at most $s^2/(2t)-s/2$ edges. Indeed, the probability that
there is a subset of size $s$ which violates this assertion is
at most
\begin{eqnarray*}
\mathbb{P} &\leq& {n \choose s}{s^2/2 \choose
s^2/(2t)-s/2}p^{s^2/(2t)-s/2} \leq \left(\frac{en}{s}\right)^s
\left(e\,\frac{s}{s-t}\,tp\right)^{s^2/(2t)-s/2} \\
 &\leq&
\left(\frac{en}{s}\right)^s
\left(\frac{1}{2e^2}\Big(\frac{n}{s}\Big)^{-2t/(s-t)}\right)^{s^2/(2t)-s/2}
\leq 2^{-s/2}=o(1),
\end{eqnarray*}
(where the $o(1)$-term tends to zero as $s$ tends to infinity).
This implies that with high probability
every subgraph of $G$ on $s$ vertices has
average degree $d \leq s/t-1$. Therefore by Tur\'an's theorem it
contains an independent set of size at least $\frac{s}{d+1} \geq t$. On the
other hand, it is well known (see, e.g., \cite{Bol}), that almost
surely the independence number of $G_{n,p}$ is bounded by
$$\alpha(G_{n,p}) \leq O\left(p^{-1}\log np\right) \leq
O\left(t\Big(\frac{n}{s}\Big)^{2t/(s-t)}\log (n/t) \right).$$ This
implies that a.s. $G_{n,p}$ satisfies the assertion of the theorem
and completes the proof. \hfill $\Box$

For the proof of Theorem \ref{upper-bound2} we
need the following lemma.
%, whose proof we include
%here for the sake of completeness.

\begin{lemma}
\label{independence} Let $G=G_{s,p}$ be a random graph, assume
$sp \rightarrow \infty$ and fix $\epsilon>0$. Then the probability
that the independence number of $G$ is at most $
\frac{\epsilon}{p}\log (sp)$ is less than $e^{-s
(sp)^{1-3\epsilon/2}}$.
\end{lemma}

\noindent {\bf Proof (of lemma).}\, 
Let $k=\frac{\epsilon}{p}\log (sp)$. To
prove the lemma we use the standard greedy algorithm which
constructs an independent set by examining the vertices of the graph
in some fixed order and by adding a vertex to the current independent
set whenever possible. The behavior of this algorithm for random
graphs can be analyzed rather accurately (see, e.g., \cite{Bol}).
At iteration $i$ of our procedure we use the greedy algorithm to find
a maximal (with respect to inclusion) 
independent set $I_i$ in the remaining vertices of $G$.
If $|I_i| \geq k$ we stop. Otherwise we delete the vertices of $I_i$
from $G$ and continue. We stop when the number of remaining
vertices drops below $s/2$. 
Note that during iteration number $i$ we only expose edges incident to 
$I_i$, therefore the remaining vertices still form a truly random graph.
Given a set $I$ the probability that a fixed vertex of $G$ 
is adjacent to some vertex of
$I$ is $1-(1-p)^{|I|}$. Therefore the probability that a fixed 
$I$ is maximal, where $|I| \leq k $, is 
at most $\big(1-(1-p)^k\big)^{s/2}$. 
By definition, when the iteration fails, the random graph must contain a
maximal independent set of size less than $k$ (and hence also a set of size exactly
$k$ so that every remaining vertex is adjacent to at least one of its members).
Thus the probability of this event is at most ${s \choose k}\big(1-(1-p)^k\big)^{s/2}$. 
Moreover, the
outcomes of different iterations depend on disjoint sets of edges
and therefore are independent.  Finally note that if $G$ has no
independent set of size $k$, then the number of iterations is at
least $s/(2k)$. This implies that the probability of such an event is
bounded by
\begin{eqnarray*}
\mathbb{P} &\leq& \left({s \choose k}\big(1-(1-p)^k\big)^{s/2}\right)^{s/(2k)}
\leq \left(\frac{es}{k}\right)^{s/2}e^{-\frac{s^2}{4k}(1-p)^k}\\
\hspace{4.1cm} &\leq& e^{-\Omega\big(s
\frac{(sp)^{1-\epsilon}}{\log sp}\big)+s\log (sp)} \leq e^{-s
(sp)^{1-3\epsilon/2}}\,. \hspace{5.2cm}\Box
\end{eqnarray*}

\noindent {\bf Proof of Theorem \ref{upper-bound2}.}\, 
Suppose
that $(s/t)^{1-\delta} \geq \log n$ for some fixed $\delta >0$ and
consider the random graph $G=G_{n,p}$ with $p=\frac{\delta \log
(s/t)}{2t}$. (Note that $p<1$ as $s \leq e^{2t}$.)
As already mentioned above a.s. the independence
number of this graph is bounded by $O\big(\frac{1}{p}\log
(np)\big)=O\big(\frac{t}{\log (s/t)} \log (n/t)
\big)$. Here we used that $\log (s/t)<s/t<n/t$ and hence 
$\log\big(\frac{n}{t} \log (s/t)\big)<2 \log (n/t)$. Also,
by Lemma \ref{independence} (with
$\epsilon=\delta/2$), the probability that $G$ contains an induced
subgraph of order $s$ with no independent set of size
$\frac{\delta}{2p}\log (sp) \geq t$ is at most
$${n \choose s} e^{-s (sp)^{1-3\delta/4}} \leq n^s e^{-s  \log^{1+\delta/4} n}=o(1) .$$
Here we used that $\log(s/t) \rightarrow \infty$ and hence
$(sp)^{1-3\delta/4} \geq (s/t)^{1-3\delta/4} \geq
\log^{1+\delta/4}n$. Therefore with high probability $G$ satisfies
the first assertion of the theorem.

To prove the second part of the theorem suppose that
$s/t=\Omega(\log n)$ and $\log t \geq \log^{\gamma} n$ for some
constant $\gamma >0$. By the previous paragraph we can also assume
that $s/t \leq \log^2 n$. Let $G=G_{n,p}$ with $p=\frac{\gamma
\log (s/t)}{640t}$. Again we have that a.s. $\alpha(G) \leq
O\big(\frac{t}{\log (s/t)} \log (n/t) \big)$.
Since $sp= \Omega(\log n \log \log n)$, the probability that there
is a subset of size $s$ in $G$ which spans at least $2 s^2p $
edges is bounded by
$${n \choose s} {s^2/2 \choose 2s^2p} p^{2s^2p} \leq
n^s \left(\frac{e}{4p}\right)^{2s^2p}p^{2s^2p}= n^s
\left(\frac{e}{4}\right)^{2s^2p} \leq e^{s\log n} e^{-\Omega(s
\log n \log \log n)}= o(1).$$ Also, since $s/t \leq sp \leq
t^{o(1)}$, we have that the probability that $G$ contains a subset
of order $k=\Theta(sp)$ which spans more than $k^{2-\gamma/2}$
edges is at most
\begin{eqnarray*}
{n \choose k} {k^2/2 \choose k^{2-\gamma/2}} p^{k^{2-\gamma/2}}
&\leq& n^k (k p)^{k^{2-\gamma/2}}\leq e^{k\log n}
e^{-\Omega(k(s/t)^{1-\gamma/2}\log t)}\\ &\leq& e^{k\log n}
e^{-\Omega(k\log^{1+\gamma/2} n)}=o(1).
\end{eqnarray*}

Let $H$ be an induced subgraph of $G$ of order $s$. Since the
number of edges in $H$ is a.s. at most $2s^2p$ it contains an
induced subgraph $H'$ of order $s/2$ with maximum degree 
at most $d =
8sp$. By the above discussion, we also have that the neighborhood
of every vertex of $H'$ spans at most $d^{2-\gamma/2}$ edges and
therefore through every vertex in $H'$ there are at most $t \leq
d^{2-\gamma/2}$ triangles. Now we can use the known estimate
(Lemma 12.16 in \cite{Bol}, see also \cite{AKS} for a more general
result) on the independence number of a graph containing a small number
of triangles. It implies that
\begin{eqnarray*}
\alpha(H') &\geq& 0.1 \frac{|V(H')|}{d}\left(\log
d-\frac{1}{2}\log t\right) \geq 0.1 \frac{s/2}{d}\left(\log
d-\frac{1}{2}\log d^{2-\gamma/2} \right)\\
&=&\frac{\gamma s}{80d}\log d \geq \frac{\gamma s}{640 sp}\log
(8sp) \geq t.
\end{eqnarray*}
This shows that $G$ a.s. satisfies the second assertion of our
theorem and completes the proof.
\hfill $\Box$

\section{Remarks}

Theorems \ref{lower-bound1} and \ref{upper-bound1} show that if 
$s/t=O(1)$ and $t>1$ then $f(n,s,t)=n^{\Theta(1)}$, whereas if $s/t \gg 1$
and $t=n^{o(1)}$ 
then $f(n,s,t) =n^{o(1)}$, and if $s/t \geq \Omega(\log n)$
and $t \leq (\log n)^{O(1)}$ 
then $f(n,s,t) \leq (\log n)^{O(1)}.$

Theorems \ref{lower-bound2} and \ref{upper-bound2} determine the asymptotic behavior
of the function $f(n,s,t)$ up to a constant factor for a wide range of the parameters.
We do not specify here all this range, and only observe that 
in particular,  for every
fixed $\mu>0$, $f(n,\log^{2+\mu} n, \log n)=\Theta(\frac{\log^2 n}{\log \log n}).$
For $\mu=1$ this implies (\ref{e12}). The estimate (\ref{e11})
follows from Theorems \ref{lower-bound2} and \ref{upper-bound1}.

The connection between the Ramsey-type question described in Section 1 and the function
$f$ is the following simple fact.

\noindent
{\bf Claim:}\, If 
\begin{equation}
\label{e41}
n/2 >s>t ~~\mbox{ and }~~ (t-1)f(n/2,s,t) \geq s, 
\end{equation}
then there is no
graph on $n$ vertices in which every induced subgraph on $s$ vertices contains
a clique of size $t$ and an independent set of size $t$.

\noindent
{\bf Proof:}\, Assuming there is such a graph $G$, observe that by the definition of
$f$ it contains an independent set $I_1$ of size  $f=f(n/2,s,t)$. Omit this set,
and observe that the induced graph on the remaining vertices, assuming there are at
least $n/2$ of them, contains another independent set $I_2$ of size $f$.
Repeating this process $t-1$ times (or until the union of the 
independent sets obtained is of
size at least $n/2>s$), we get an induced subgraph of $G$ on 
$min\{ n/2,(t-1)f \} \geq s$ vertices,
which is $(t-1)$-colorable (as it is the union of $t-1$ independent sets), and hence
cannot contain a clique of size $t$. This is a contradiction, proving the assertion of
the claim.

In particular, for $t=\log n$ and $s=c \log^3 n/\log \log n$, where $c$ is a
sufficiently small positive absolute 
constant, it is not difficult to check that the assumption
in (\ref{e41}) holds, by Theorem \ref{lower-bound2}.

It will be interesting to close the gap between our upper and lower bounds 
for the function $f(n,s,t)$. It will also be interesting to know more about
the Ramsey-type question of Erd\H{o}s and Hajnal described in the introduction,
and in particular, to decide if there exists a graph on $n$ vertices in which
every induced subgraph on, say, $\log^{100} n$ vertices, contains a clique of size at
least $\log n$ and an independent set of size at least $\log n$.

\end{document}